\documentclass[11pt,leqno,twoside]{article}
\usepackage{latexsym,theorem}
\input latex.sty
\usepackage{amssymb,amsmath}    

\newtheorem{Exercise}{Exercise}[section]
\newtheorem{Rule}{Rule}[section]

\def\d{\mathord{\rm d}}
\def\E{\mathord{\rm E}}
\def\P{\mathord{\rm P}}
\def\trace{\mathop{\rm trace}\nolimits}

\def\bra#1{\left \langle #1 \right |}
\def\ket#1{\left | #1 \right \rangle}
\def\braket#1#2{\left \langle #1 \mid  #2 \right \rangle}
\def\braopket#1#2#3{\left \langle #1 \mid  #2 \mid  #3\right \rangle}

\begin{document}

\Title
{Asymptotics in Quantum Statistics}
 
\shorttitle
{Quantum Asymptotics}

\Author
{Richard~D.~Gill}

\affiliation
{Mathematical Institute, University Utrecht, and\\
Eurandom, Eindhoven.}

\Abstract
{Observations or measurements taken of a quantum system (a small number 
of fundamental particles) are inherently random. If the state of the 
system depends on unknown parameters, then the distribution of the outcome
depends on these parameters too, and statistical inference problems 
result. Often one has a choice of what measurement to take,
corresponding to different experimental set-ups or settings of measurement
apparatus. This leads to a design problem---which measurement is best for
a given statistical problem. This paper gives an introduction to this
field in the most simple of settings, that of estimating the state of a
spin-half particle given $n$ independent copies of the particle.
We show how in some cases asymptotically optimal measurements can be
constructed. Other cases present interesting open problems, connected
to the fact that for some models, quantum Fisher information is
in some sense non-additive. In physical terms, we have non-locality
without entanglement.}

\subject
{62F12, 62P35.}


\section{Introduction}
\label{s1:intro}
The fields of quantum statistics and quantum probability have a reputation for being 
esoteric. However, in our opinion, quantum mechanics is a fascinating source of probabilistic and statistical
models, unjustly little known to `ordinary' statisticians and probabilists.

Quantum mechanics has two main ingredients: one deterministic, one random.
In isolation from the outside world a quantum system evolves deterministically
according to Schr\"odinger's equation. That is to say, it is described by
a state or wave-function whose time evolution is the (reversible) solution of a 
differential equation. On the other hand when this system comes into
interaction with the outside world, as when for instance measurements are
made of it (photons are counted by a photo-detector, tracks of particles
observed in a cloud chamber, etc.\@) something random and irreversible takes place. The state of the system makes a random jump and the outside world
in some way contains a record of the jump.
 From the state of the system at the time of the interaction one can read
off, according to certain rules, the probability distribution of the 
macroscopic outcomes
and the new state of the system. (See Penrose, 1994, for an eloquent
discussion of why there is something paradoxical in the peaceful coexistence
of these two principles; and see Percival (1998) for interesting stochastic 
modifications to Schr\"odinger's equation which might offer some 
reconciliation).\footnote{Also highly recommended: Sheldon Goldstein,
`Quantum mechanics without observers', Physics Today, March, April 1998;
letters to the editor, Physics Today, February 1999.}

Till recently most predictions made from quantum theory involved such
large numbers of particles that the law of large numbers takes over and
predictions are deterministic. However technology is rapidly
advancing to the situation that really small quantum systems can be
manipulated and measured (e.g., a single ion in a vacuum-chamber, or
a small number of photons transmitted through an optical communication
system). Then the outcomes definitely are random. The fields of
quantum computing, quantum communication, and quantum cryptography 
are rapidly developing and depend on the ability to manipulate really
small quantum systems.
Theory and conjecture are much further than experiment and technology,
but the latter are following steadily.

In this paper we will introduce as simply as possible the model of quantum
statistics and consider the problem of how best to measure the state of
an unknown spin-half system. We will survey some recent results, in
particular, from joint work with O.E.~Barndorff-Nielsen and with 
S.~Massar (Barndorff-Nielsen and Gill, 1998; Gill and Massar, 1998).
This work has been concerned with the problem, posed by Peres and Wootters
(1991): can more information be obtained about the common state of $n$ 
identical
quantum systems from a single measurement on the joint system formed by
bringing the $n$ systems together, or does it suffice to combine separate
measurements on the separate systems? A useful tool for our studies is
the quantum Cram\'er-Rao bound with its companion notion of quantum 
information, introduced by C.W.~Helstrom in a sequence of papers in the
sixties and later refined by among others A.S.~Holevo. 

Quantum statistics mainly consists of exact results in various rather special 
models, see the books of Helstrom (1976) and Holevo (1982). Just as
in ordinary statistics, the Cram\'er-Rao bound on the variance of an
unbiased estimator is rarely achieved
exactly (only in so-called quantum exponential models). In any case, one 
would not
want in practice to restrict attention to unbiased estimators only.
There are results on optimal invariant methods, but again, not many models
have the structure that these results are applicable and even then the
restriction to invariant statistical methods is not entirely compelling.

One might hope that asymptotically it would be possible to
achieve the Cram\'er-Rao bound. However asymptotic theory is so far 
very little developed in the theory of quantum statistics, one reason being 
that the powerful modern tools of asymptotic statistics (contiguity, local
asymptotic normality, and so on) are just not available\footnote{though 
R. Rebolledo is working on a notion of quantum contiguity}
since even if we
are considering measurements of $n$ identical quantum systems, there
is no a priori reason to suppose that a particular sequence of measurements
on $n$ quantum systems together will satisfy these conditions. Here,
we make a little progress through use of the van Trees inequality
(see Gill and Levit, 1995), a Bayesian Cram\'er-Rao bound, which will
allow us to make asymptotic optimality statements without assuming or proving
local asymptotic normality. Another useful ingredient will be the recent
derivation of the quantum Cram\'er-Rao bound by Braunstein and Caves (1994),
linking quantum information to classical expected Fisher information in a particularly
neat way.

We will show that for certain problems, a new Cram\'er-Rao type 
inequality of Gill and Massar (1998)
does provide an asymptotically achievable bound to the
quality of an estimator of unknown parameters. For some other problems
the issue remains largely open and we identify situations where 
Peres and Wootter's question has an affirmative answer: there can be
appreciably more information in a joint measurement of several particles than
in combining separate measurements on separate particles. This
clarifies an earlier affirmative answer of Massar and Popescu (1995),
which turned out only for small samples to improve on separate measurements.
It also clarifies the recent findings of Vidal et al.\ (1998).

Helstrom wrote in the epilogue to his (1976) book: ``{\em Mathematical statisticians
are concerned with asymptotic properties of estimators. When the parameters
of a quantum density operator are estimated on the basis of many
independent observations, how does the accuracy of the estimates depend
on the number of the observations as that number grows very large? Under what
conditions have the estimators asymptotically normal distributions?
Problems such as these, and still others that doubtless will occur to
physicists and mathematicians, remain to be solved within the framework of
the quantum mechanical theory.}'' More than twenty years later this programme
is still hardly touched (some of the few contributions are by Brody and Hughston
(1998) and earlier papers,  and Holevo (1983))
but we feel we have made a start here.

In $20\pm\epsilon$ pages (even when $\pm\epsilon=+10$) it is difficult to give a complete introduction to the topic,
as well as a clear picture of recent results. The classic books by 
Helstrom and Holevo mentioned above are still the only books on quantum
statistics and they are very difficult indeed to read for a beginner.
A useful resource is the survey paper by Malley and Hornstein (1993).
However the latter authors, as many others, take the stance that the
randomness occuring in quantum physics cannot be caught in a standard
Kolmogorovian framework. We argue elsewhere (Gill, 1998), in a critique of
an otherwise excellent introduction to the related field of
quantum probability (K\"ummerer and Maassen, 1998), that this is
nonsense.
With more space at our disposal we would have included extensive
worked examples; however they have been replaced by {\it exercises\/} 
so that the reader can supply some of the extra pages
(but---unless you are Willem van Zwet---leave the starred exercises 
for later).

Some references which we found specially useful in getting to grasps
with the mathematical modelling of quantum phenomena are the books
by Peres (1995), and Isham (1995). To get into quantum probability,
we recommend Biane (1995) or Meyer (1986).

This introductory section continues with three subsections summarizing
the basic theory: first the mathematical model of states and measurements;
secondly the basic facts about the most simple model, namely of a 
two-state system; and thirdly the basic quantum Cram\'er-Rao bound.
That third subsection finishes with a glimpse of how one might do
asymptotically optimal estimation in one-parameter models: in a preliminary
stage obtain a rough estimate of the parameter from a small number of
our $n$ particles. Estimate the so-called
quantum score at this point, and then go on to measure it in the second
stage on the remaining particles.
Section 2 states a recent new version of the quantum Cram\'er-Rao bound
which makes precise how one might trade information between different
components of a parameter vector. Section 3 outlines the procedure for
asymptotically optimal estimation of more than one parameter, again a
two-stage procedure. This is work `in progress', so some results are
conjectural, imprecise, or improvable. In a final short section we
try to explain how some of our results are connected to the strange
phenomenon of {\it non-locality without entanglement\/}, a hot
topic in the theory of quantum information and computation.

\subsection{The basic set-up}
\label{ss1.1:set-up}
Quantum statistics has two basic building blocks:
the mathematical specification of the state of a quantum system,
to be denoted by $\rho=\rho(\theta)$ as it possibly depends on
an unknown parameter $\theta$,
and the mathematical specification of the measurement, denoted 
by $M$, to be carried out on that
system. We will give the recipe for the probability
distribution of the observable outcome 
(a value $x$ of a random variable $X$ say)
when measurement $M$ is carried out on a system in state $\rho$. Since the
state $\rho$ depends on an unknown parameter $\theta$, the distribution of
$X$ depends on $\theta$ too, thereby setting a statistical problem of
how best to estimate or test the value of $\theta$. Since we may in practice
have a choice of which measurement $M$ to take, we have a design problem of
choosing the best measurement for our purposes. (There is also a recipe
for the state of the system after measurement, depending on the outcome,
but we do not need it here; Bennett et al., 1998).

For simplicity we restrict attention to finite-dimensional quantum systems.
The state of a $d$-dimensional quantum system will be modelled or specified
by a $d\times d$ complex matrix $\rho$ called the density matrix of the system.
For instance, when we measure the spin of an electron in a particular direction
only two different values can occur, conventionally called `up' and `down'.
One could call this a two-state system, we need a $d=$two-dimensional state space.
Similarly if we measure whether a photon is polarized in a particular direction
by passing it through a polarization filter, it either passes or does not
pass the filter. Again, polarization measurements on a single photon
can be discussed in terms of a two-dimensional system. If we 
consider the spins of $n$ electrons, then $2^n$ different outcomes are 
possible and the system of $n$ electrons together (or rather, their spins),
is described by a $d\times d$ matrix $\rho$ with $d=2^n$.

\begin{Definition}[Density matrix]
\label{d:density matrix}
The density matrix $\rho$ of a $d$-di\-men\-sion\-al
quantum system is a $d\times d$ self-adjoint, nonnegative matrix of trace 1.
\end{Definition}
`Self adjoint' means that $\rho^*=\rho$ where the $*$ denotes the complex
conjugate and transpose of the matrix. That $\rho$ is nonnegative means that
$\psi^*\rho\psi\ge0$ for all column vectors $\psi$ (since $\rho$ is self-adjoint
this quadratic form is a real number). We often use the Dirac bra-ket notation
whereby $\ket\psi$ (called a ket) is written for the column vector $\psi$ and
$\bra\psi$ (a bra)
is written for its adjoint, the row vector containing the complex conjugates of
its elements. The quadratic form $\psi^*\rho\psi$ is then denoted
$\braopket\psi\rho\psi$.

It follows that the diagonal elements of a density matrix are nonnegative reals
adding up to one. Moreover by the eigenvalue-eigenvector decomposition of
self-adjoint matrices we can write $\rho=\sum_i p_i\ket i \bra i$
where the kets $\ket i$ are the orthonormal eigenvectors of $\rho$,
$\braket i j =\delta_{ij}$, and the $p_i$ are the eigenvalues: nonnegative
real numbers adding up to one. One says that the density matrix $\rho$
represents the {\it mixed\/} state obtained by taking with probability $p_i$
the system in the {\it pure\/} state $\ket i$. The state vector of a pure
state is also called a wave-function.

\begin{Definition}[Measurement]
\label{d:measurement}
A measurement $M$ on a $d$-dimensional quantum
system taking values $x$ in a measurable space $(\mathcal X,\mathcal A)$ is
specified by an {\em operator-valued probability measure\/} 
or {\em oprom\/} for short, 
that is, a collection
of self-adjoint matrices $M(A):A\in\mathcal A$ such that
\begin{enumerate}
\item $M(\mathcal X)=\mathbf 1$, the identity matrix,
\item Each $M(A)$ is non-negative.
\item For disjoint $A_i$, $M(\cup_i A_i)=\sum_i M(A_i)$.
\end{enumerate}
\end{Definition}
Note that these three rules are the ordinary axioms of a probability measure
on $(\mathcal X,\mathcal A)$, except that the measure takes values in the
self-adjoint matrices instead of the real numbers.
The sample space $\mathcal X$ might be the
real numbers or a subset thereof, with the Borel sigma algebra, 
but it could also be anything else.

Measurements are often called {\em generalised measurements\/}, to contrast them
with a special subclass of measurements called {\em simple measurements\/}
which we will introduce in a moment. In the literature the abbreviations
`povm' (positive operator valued measure) and `pom' (probability operator
matrices) are often used, which we however find inaccurate.

Now we can give the so-called {\em trace-rule\/} telling us the probability distribution of
the random outcome $X$ when $M$ is used to measure $\rho$:

\begin{Definition}[trace rule]
\label{d:trace rule}
The probability distribution of the outcome $X$ is given by
\begin{equation}
\label{e:trace rule}
\Pr\{X\in A\}~=~\trace(\rho M(A)),\qquad A\in\mathcal A
\end{equation}
\end{Definition}

\begin{Exercise}[legitimacy of trace rule]
Prove that \eqref{e:trace rule} indeeds defines a probability measure on
$\mathcal X,\mathcal A$.
\end{Exercise}

One can argue from basic principles of quantum mechanics that however one 
measures a quantum system, the result must be an affine mapping from density 
matrices to the space of probability distributions
on the outcome space. It is a theorem that any such mapping can be
represented by an oprom. Thus the class of oproms contains all
conceivable measurements. 
On the other hand, as we will see later, any oprom can be 
realised by some concrete experimental set-up, at least in principle,
so the definition captures exactly what it should.

A special kind of measurements plays a key role in theory and practice,
these are the so-called {\em simple measurements\/} defined as follows:
\begin{Definition}[Simple measurement]
\label{d:simple measurement}
A simple measurement $\Pi$ on a $d$-dimensional quantum
system taking values $x$ in a measurable space $(\mathcal X,\mathcal A)$ is
a measurement such that each $\Pi(A)$ is idempotent, i.e., is a projector
onto a subspace of $\mathbb C^d$.
\end{Definition}

It follows that the measurement takes on
at most $d$ different values, i.e., there exist $x_1,\dots,x_k\in\mathcal X$
with $k\le d$ such that $\Pi(\{x_1,\dots,x_k\})=\mathbf 1$. Writing $\Pi(x_i)$
as abbreviation for $\Pi(\{x_i\})$ the matrices $\Pi(x_i)$ project onto $k$
orthogonal subspaces of $\mathbb C^d$ together spanning the whole space.
Let us now define a self-adjoint matrix $X$ (not to be confused with
the random variable $X$ representing the outcome of the measurement)
by $X=\sum_i x_i\Pi(x_i)$. Then the $x_i$
are the eigenvalues of $X$ and the $\Pi(x_i)$ project onto the eigenspaces.
Conversely, given a self-adjoint matrix $X$ one can construct a corresponding
simple measurement or projector-valued probability measure. In this role we
call $X$ an {\em observable\/}. It follows that the expected value of
the outcome of a measurement of $X$ is given by $\trace(\rho X)$. For an
ordinary real function $f$ (e.g., square, inverse, logarithm,\dots) one defines the same function
of the observable $X$ by $f(X)=\sum_i f(x_i)\Pi(x_i)$,
and the expected value of the outcome of a measurement of the observable 
$f(X)$ is $\trace(\rho f(X))$.

Simple measurements are often called von Neumann measurements. We will 
occasionally use the term `proprom' (projector-valued probability measure).
Physicists generally agree that any simple measurement could in principle be
implemented in practice. 

`Between measurements' a quantum system evolves deterministically according
to the famous Schr\"odinger equation, a differential equation for the component
pure states $\ket i$ of a given mixed system. One thinks of a measurement as
taking place instantaneously. After the measurement, the quantum system jumps
to a new state (depending on the outcome $x$); this is called `the collapse of
the wave function'. Again some simple rules specify what happens, but we will
not give them here.

If we bring two separate quantum systems together into some kind of interaction
then their future evolutions will be linked together. Measurements
can be made on the `joint system', including all the separate measurements
on each of the separate systems but many more besides.  Mathematically this
is modelled as follows:

\begin{Definition}[product system]Consider two quantum systems,
of dimension $d$ and $d'$, in states
$\rho$ and $\rho'$ respectively. Together the two are in the state
$\rho\otimes\rho'$ in 
$\mathbb C^d \otimes \mathbb C^{d'}=\mathbb C^{d\times
d'}$ where $\otimes$ denotes the tensor product (of matrices, vectors, or
spaces as appropriate).
\end{Definition}

For the reader who is not familiar with tensor products,
the tensor product of $\mathbb C^d$ with $\mathbb C^{d'}$ has as basis
the tensor product of each element of a basis of $\mathbb C^d$ with
each element of a basis of $\mathbb C^{d'}$. One can take linear combinations
of tensor products $\psi\otimes\psi'$ by expanding bilinearly in chosen 
basis' of the two spaces. Tensor products of matrices are defined in the 
natural way by how they operate on products of vectors:
$X\otimes X' \, \psi\otimes\psi'=X\psi\, \otimes\, X'\psi'$. The trace of
a tensor product of two matrices is a product of the traces.

Suppose $M$ and $M'$ are measurements on two separate quantum systems $\rho$
and $\rho'$. Then we can define a joint measurement $M\otimes M'$ on the
combined system in the obvious way, taking values in the product of the
outcome spaces of $M$ and $M'$.
\begin{Exercise}[product measurement]
\label{e:product measurement}
Show that the outcome of measurement of $M\otimes M'$ on a system in state
$\rho\otimes\rho'$ is distributed as independent realisations of measurement
of $M$ and $M'$ on $\rho$ and $\rho'$ respectively.
\end{Exercise}
However the important point is that bringing two quantum systems together
allows many more measurements than just product measurements (which as we
saw from exercise \ref{e:product measurement} are not very interesting).

Product systems are important for two main reasons. Firstly, one of the
main themes of this paper is going to be: if we have $n$ independent
systems each in the same state $\rho(\theta)$ (i.e., in identical states
all depending on the {\em same\/} unknown parameter $\theta$), can we learn 
more about $\theta$ from a joint measurement on the $d^n$ dimensional 
combined system  $\rho^{\otimes n}(\theta)$? In the next section we will
discuss some of the history and other background to this question, which has
been the subject of a series of papers in recent years.
Secondly, product systems play a role in the realisation of generalised
measurements. It is a theorem (due to Naimark) that {\em any generalised
measurement whatever can be realised by a simple measurement after a
`quantum randomisation'}. That is to say, given any measurement $M$ there
exists a so-called ancillary system in state $\rho'$ and a simple
measurement $\Pi$ on the joint system $\rho\otimes\rho'$ such that
$\trace(\rho M(A))=\trace(\rho\otimes\rho'\, \Pi(A))$ for all $A$ and
whatever $\rho$.

\subsection{Spin half}
\label{ss1.2:spin half}
In order to make the above rather abstract concepts a little more concrete,
let us go to the most simple special case, $d=2$. This is the appropriate
set-up for studying spin-half systems like the electron. We will see that
we can associate the state of a spin-half system with a real vector $\vec a$
of length less than or equal to $1$ in ordinary
three dimensional space, and a simple measurement---which can take on
at most two different values---with a direction in space, or a unit vector 
$\vec u$. The trace rule
\eqref{e:trace rule} will reduce to a very simple formula involving $\vec a$ and
$\vec u$. The model applies to the famous Stern-Gerlach experiment, 
featuring in many introductory textbooks on quantum physics.
In that experiment silver
atoms were made to pass through a strongly varying magnetic field,
having a certain direction. Each
atom was either deflected upwards or downwards with respect to the direction
of the field. The deflection is due to
the spin of the outermost electron in the silver atom, which can be
characterized by a vector $\vec a$. The orientation of the
magnet determines which measurement is being taken, i.e., the value of $\vec u$.

First we take some time to study some special features of the $2\times 2$ self-adjoint 
matrices. The properties we find will greatly simplify calculations.
Let $\mathbf 1$ denote the identity matrix and define the {\em Pauli spin matrices\/} as follows:

\begin{Definition}[Pauli spin matrices]
\begin{equation}\label{e:pauli def}
\sigma_x=
\left(
\begin{array}{cc}
0 & 1 \\
1 & 0 
\end{array}
\right),\quad
\sigma_y=
\left(
\begin{array}{cc}
0 & -i \\
i & 0 
\end{array}
\right),\quad
\sigma_z=
\left(
\begin{array}{cc}
1 & 0 \\
0 & -1 
\end{array}
\right).
\end{equation}
\end{Definition}
These three matrices are self adjoint, each have trace zero and determinant
minus one, hence have eigenvalues $\pm 1$. They satisfy (check this yourself!)
\begin{equation}\label{e:pauli props}
\begin{array}{rcl}
\sigma_x\sigma_y=-\sigma_y\sigma_x&=&i\sigma_z,\\
\sigma_y\sigma_z=-\sigma_z\sigma_y&=&i\sigma_x,\\
\sigma_z\sigma_x=-\sigma_x\sigma_z&=&i\sigma_y,\\
\sigma_x^2=\sigma_y^2=\sigma_z^2&=&{\mathbf 1}.
\end{array}
\end{equation}
An arbitrary self-adjoint $2\times 2$ complex matrix has to be of the form
\begin{equation}
X=\left(
\begin{array}{cc}
u+z & x-iy \\
x+iy & u-z 
\end{array}
\right)
\end{equation}
where $x,y,z,u$ are uniquely determined real numbers. 
Thus we can write
\begin{equation}
X=u{\mathbf 1}+x \sigma_x + y \sigma_y + z \sigma_z.
\end{equation}

Specializing to density matrices, the requirement that $\trace\rho=1$
imposes the condition that $u=\frac12$. The requirement that $\rho$ is 
nonnegative is equivalent to its determinant being nonnegative, or
$u^2-z^2-x^2-y^2\ge 0$, or $x^2+y^2+z^2\le{\frac 1 2}^2$. It is convenient
to write
\begin{equation}
\rho=\rho(\vec a)=\frac12({\mathbf 1}+\vec a\cdot\vec\sigma)
\end{equation}
where $\vec a=(a_x,a_y,a_z)\in \mathbb{R}^3$ and satisfies
\begin{equation}
\|\vec a\|^2=a_x^2+a_y^2+a_z^2\le 1
\end{equation}
while $\vec\sigma=(\sigma_x,\sigma_y,\sigma_z)$ (a vector of matrices)
and `$\cdot$' denotes the inner-product. Thus the space of density matrices 
of a two-dimensional quantum system can be represented by the closed unit ball 
$B$ in
three dimensional Euclidean space. The sphere $S$, or surface of the unit ball,
corresponds to density matrices $\frac12({\mathbf 1}+\vec a\cdot\vec\sigma)$
with $\|\vec a\|^2= 1$
which are singular since their determinant is
zero. Such a density matrix has therefore eigenvalues $0$ and $1$.
It represents a so-called {\em pure state\/}.

The density matrix of a pure state is a projector matrix, projecting 
onto a one-dimensional subspace
of $\mathbb{C}^2$. Letting $\vec u$ denote a unit vector in $\mathbb{R}^3$,
let us write $\Pi(\vec u)=\rho(\vec u)=\frac12({\mathbf 1}+\vec u\cdot\vec\sigma)$ for this matrix. Check 
using \eqref{e:pauli props} that
 $\Pi(\vec u)$ is idempotent, and that $\Pi(\vec u)$
and $\Pi(-\vec u)$ commute (in fact, their product is the zero matrix)
and add to the identity matrix! Thus the projectors  $\Pi(\vec u)$
and $\Pi(-\vec u)$ project onto two orthogonal one-dimensional subspaces
of $\mathbb{C}^2$. We will specify these spaces exactly in a moment.
The only other projector matrices are $\mathbf 0$ and $\mathbf 1$, projecting
onto the trivial subspace and the whole space of $\mathbb C^2$ respectively.

It follows that for an arbitrary density matrix $\rho=\rho(\vec a)$,
defining the unit vector $\vec u=\vec a/\|\vec a\|$ 
and the probabilities $\alpha=\|\vec a\|$, $\beta=1-\alpha$, we have
\begin{equation}
\begin{array}{rcl}
\rho(\vec a)=\frac12({\mathbf 1}+\vec a\cdot\vec \sigma)
&=& \|\vec a\| \rho(\vec a/\|\vec a\|)
+(1-\|\vec a\|) \rho(-\vec a/\|\vec a\|)\\
&=&\alpha \rho(\vec u)
+\beta \rho(-\vec u).
\end{array}
\end{equation}
It has eigenvalues $\alpha$ and $\beta$, and its eigenvectors,
column vectors in $\mathbb{C}^2$,
generate the spaces onto which  $\Pi(\vec u)$ and
$\Pi(-\vec u)$ project. One may consider the state $\rho(\vec a)$ as the
mixture, with probabilities $\alpha$ and $\beta$, of the pure states
$\rho(\vec u)$ and $\rho(-\vec u)$ (though this is only one of
many representations of $\rho$ as a mixture of pure states).

So what are these spaces exactly? The vector $\vec u$ is a point on the
unit sphere in $\mathbb{R}^3$. Let $\theta$ and $\phi$ denote its polar
coordinates, where $\theta\in[0,\pi]$ is the latitude measured from the 
North pole ($z$-axis) and $\phi\in[0,2\pi)$ is the longitude, measured from 
the $x$-axis. (We should really say co-latitude rather than latitude).
Thus $\vec u=(\sin\theta\cos\phi,\sin\theta\sin\phi,\cos\theta)$.
Define the column vector $\ket\psi=\ket{\psi(\theta,\phi)}$ in $\mathbb{C}^2$ by
\begin{equation}
\ket{\psi(\theta,\phi)}=\left(
\begin{array}{c}
e^{-i\phi/2}\cos(\theta/2)\\
e^{i\phi/2}\sin(\theta/2)
\end{array}\right).
\end{equation}
Note that $\braket\psi\psi=1$ while
\begin{equation}
\begin{array}{rcl}
\ket\psi\bra\psi&=&
\left(
\begin{array}{c}
e^{-i\phi/2}\cos(\theta/2)\\
e^{i\phi/2}\sin(\theta/2)
\end{array}\right)
\left(
\begin{array}{cc}
e^{i\phi/2}\cos(\theta/2)
e^{-i\phi/2}\sin(\theta/2)
\end{array}\right)
\\
&=&
\left(\begin{array}{cc}
\cos^2(\theta/2) & e^{-i\phi}\cos(\theta/2)\sin(\theta/2)
\\
e^{i\phi}\cos(\theta/2)\sin(\theta/2) & \sin^2(\theta/2)
\end{array}
\right)
\\
&=&
{\frac 1 2}
\left(
\begin{array}{cc} 
1+\cos(\theta) & (\cos\phi-i\sin\phi) \sin\theta \\
(\cos\phi+i\sin\phi) \sin\theta & 1-\cos\theta
\end{array}
\right)
\\
&=&{\frac 1 2}({\mathbf 1}+\vec u\cdot\vec\sigma)=\Pi(\vec u).
\end{array}
\end{equation}
{\it Any\/} complex vector $\ket\xi$ of length $1$ can be written
as $e^{i\alpha}\psi(\theta,\phi)$ for some $\alpha\in[0,2\pi)$ and polar 
coordinates $\theta,\phi$. Note that $\ket\xi\bra\xi=\ket\psi\bra\psi=\Pi(\vec u)$, and that $\ket{\psi(\theta,\phi)}$ and $\ket{\psi(\pi-\theta,\phi+\pi)}$
are orthogonal. The corresponding points on the unit sphere 
are opposite to one another. Combining these facts we obtain:

\begin{Rule}[Spin-half density matrices, projectors]
The density matrix $\rho(\vec a)$,
where $\vec a$ is a point in the unit ball in $\mathbb{R}^3$, has eigenvalues
$\|\vec a\|$ and $1-\|\vec a\|$ and normalized
eigenvectors $\ket{\psi(\theta,\phi)}$,
$\ket{\psi(\pi-\theta,\phi+\pi)}$, where $\theta$ and $\phi$ are the polar
coordinates of $\vec u=\vec a/\|\vec a\|$. The projector matrix $\Pi(\vec u)$
projects onto the one-dimensional subspace of $\mathbb{C}^2$ spanned by 
$\ket{\psi(\theta,\phi)}$. The projector onto the space orthogonal to this,
spanned by $\ket{\psi(\pi-\theta,\phi+\pi)}$, is $\Pi(-\vec u)$.
\end{Rule}

Let $\vec u$ and $\vec v$ be two unit vectors in $\mathbb{R}^3$ and write $\ket{\vec u}$ and
$\ket{\vec v}$ for the corresponding unit vectors in $\mathbb{C}^2$; 
so $\ket{\vec u}$ is an 
abbreviation for $\ket{\psi(\theta,\phi)}$ where $\theta,\phi$ are the polar 
coordinates of $\vec u$. Since $\Pi(\vec u)=\ket{\vec u}\bra{\vec u}$
we see that $\trace\Pi(\vec u)\Pi(\vec v)=\trace\ket{\vec u}\bra{\vec u}
\ket{\vec v}\bra{\vec v}=\braket{\vec v}{\vec u}\braket{\vec u}{\vec v}=
|\braket{\vec u}{\vec v}|^2$. On the other hand, using the properties
\eqref{e:pauli props} of the Pauli matrices, one readily computes
$\trace\Pi(\vec u)\Pi(\vec v)={\frac 1 2}(1+\vec u\cdot \vec v)$. Now
$\vec u\cdot \vec v$ is the cosine of the angle between the vectors
$\vec u$ and $\vec v$, hence ${\frac 1 2}(1+\vec u\cdot \vec v)$ is the
squared cosine of half the angle between $\vec u$ and $\vec v$.

\begin{Rule}[Calculation rule]
The absolute value of the squared inner pro\-d\-uct between the
complex vectors $\ket{\vec u}$ and $\ket{\vec v}$ in $\mathbb{C}^2$ is the
squared cosine of half the angle between the corresponding unit vectors
$\vec u$ and $\vec v$ in $\mathbb{R}^3$. In particular, opposite points on the
unit sphere correspond to orthogonal vectors in $\mathbb{C}^2$.
\end{Rule}

We can now describe the probability distributions of all {\em simple\/} measurements
of a spin-half system.

The state of the system is modelled by a $2\times 2$ density 
matrix of the form $\rho(\vec a)={\frac 1 2}({\mathbf 1}+\vec a\cdot\vec
\sigma)$ where $\vec a$ is a point in the closed unit ball in $\mathbb{R}^3$.

The non-trivial simple measurements take on just two different values.
Consider a simple measurement $M=\Pi$ taking values in a set $\mathcal X$
consisting of just two elements, let's call these elements $\pm 1$.
The measurement is determined by the two projectors $\Pi(\pm)$, which 
should project onto orthogonal one-dimensional subspaces of 
$\mathbb{C}^2$. Each subspace is generated by a vector of the form
$\ket{\vec u}$ for some $\vec u$ on the unit sphere, and the associated
projectors are $\Pi(\vec u)$. Recall that opposite points $\pm\vec u$
on the unit sphere correspond to orthogonal vectors $\ket{\pm\vec u}$
in $\mathbb{C}^2$, and hence to orthogonal projectors $\Pi(\pm\vec u)$.
Thus a projector-valued probability measure for a simple
measurement with values in $\mathcal X$ is given by
$M(\pm 1)=\Pi(\pm\vec u)=\frac12({\mathbf 1}\pm\vec u\cdot\vec \sigma)$
for some $\vec u$.

We apply the trace rule \eqref{e:trace rule} to compute the probabilities
of the two outcomes $\pm 1$ when the simple measurement $M(\pm 1)=\Pi(\pm\vec u)$
is carried out on a system in the state $\rho(\vec a)=
{\frac 1 2}({\mathbf 1}+\vec a\cdot\vec\sigma)$.
Using the properties \eqref{e:pauli props} of the Pauli matrices, the reader
should verify that these probabilities are
\begin{equation}\label{e:stern gerlach probs}
\trace\rho(\vec a)\Pi(\pm\vec u)={\frac 1 2}(1\pm\vec a\cdot \vec u).
\end{equation}

Using further rules for the state of the system after measurement,
it turns out that after measurement the system is in the
pure state $\rho(\pm\vec u)$ according to the outcome $\pm 1$.
One can therefore go on to compute probabilities of the series of outcomes
of a series of simple measurements carried out on one particle.

In the Stern-Gerlach experiment, the initial state of the silver atom is described
by the density matrix $\rho(\vec 0)=\frac12{\mathbf 1}$.  One can think of this state as
corresponding to an electron having spin in a random direction $\vec u$ uniformly distributed over
the unit sphere. Indeed, if one takes the mean of $\rho(\vec u)=\frac12({\mathbf 1}+\vec u\cdot
\vec\sigma)$ with $\vec u$ uniformly distributed over the sphere, the matrix $\frac12{\mathbf 1}$
results (though this representation of the `completely random' state $\rho(\vec 0)$ as
a mixture of pure states is not unique; one also finds this state as the result
of choosing with equal probabilities $\frac12$ an electron in either of the orthogonal
pure states $\ket{\pm\vec u}$).

\begin{Exercise}[A generalised measurement of spin-half system]
Let 
\\
$M(A)=\int_A\Pi(\vec u)\d\vec u/2\pi$ where
$\d\vec u$ denotes integration with respect to Lebesgue surface measure on
$S$. Show that $M$ is a generalized measurement on a spin-half system with values
in $S$, and compute the distribution of the outcome of this measurement
on the system $\rho(\vec a)$. This measurement would be physically
realised by somehow coupling the spin-half system with a particle moving
on the sphere and measuring the position of that particle.
\end{Exercise}

\begin{Exercise}[A generalized measurement of $n$ spin-half systems*]
For the state-space $(\mathbb C^2)^{\otimes n}$
define ${\ket{\vec u}}_n=\ket{\vec u}\otimes\cdots\otimes\ket{\vec u}$
and define $\Pi_n(\vec u)={\ket{\vec u}}_n{\bra{\vec u}}_n$.
Define $M(A)=(n+1)\int_A\Pi_n(\vec u)\d\vec u/4\pi$ and
show that $M(S)$ is the projector onto the $n+1$ dimensional
subspace of vectors, invariant under permutation of the $n$ components of
$(\mathbb C^2)^{\otimes n}$. Call this subspace $\mathcal S_n$ and
note that $\trace\rho^{\otimes n} \Pi_{\mathcal S_n}=1$. Show that $M$
defines a generalized measurement on $n$ identical copies of a spin-half system 
with values in $S$, and compute the distribution of the outcome of this measurement
on the system $\rho(\vec v)$.
\end{Exercise}

A Stern-Gerlach magnet oriented in the direction $\vec u$ implements the simple measurement
$M(\pm 1)=\Pi(\pm\vec u)$. Since for $\vec a=\vec 0$ the probabilities 
\eqref{e:stern gerlach probs} both equal $\frac12$, one will find electrons with spin in
the directions $\pm\vec u$ with equal probabilities. Electrons in the emerging
`$+$' beam are in the pure state $\rho(\vec u)$. Sending them through a Stern-Gerlach device
with orientation $\vec v$ splits them again, now with probabilities 
$\frac12(1\pm\vec u
\cdot\vec v)$ (the squared cosine of half the angle between
the directions $\vec u$ and $\vec v$)
into two beams of electrons in the states $\rho(\pm\vec v)$,
 and so on.

If the electrons started out in the arbitrary mixed state $\rho(\vec a)$ then the first
Stern-Gerlach magnet splits them into two output beams in the pure states $\rho(\pm\vec u)$
in the proportions $\frac12(1\pm\vec a \cdot \vec u)$. So if $\vec a$ was unknown, we do
learn something about it from counting the numbers of electrons in each beam. Further
operations on the output beams however will not teach us any more as the state of the
electrons in either output beam no longer depends on $\vec a$.

If we are allowed to measure a large number of electrons each in the same mixed state
$\rho(\vec a)$, we see that a large number of Stern-Gerlach measurements  in three linearly 
independent directions will enable us to determine $\vec a$. The question we
will study in the rest of the paper is: what is the best way to do this?
Will it suffice to use simple measurements on separate particles or can we
do better by using more sophisticated measurements, in particular, joint
measurements on several particles simultaneously?

One can consider rotating a given coordinate system in $\mathbb R^3$
in such a way as to transform the vectors $\vec a$ and $\vec u$
representing a state or a simple measurement into convenient choices,
e.g., we will in the future claim that `without loss of generality
$\vec a = (0,0,a_3)$' which makes $\rho(\vec a)$ a diagonal matrix.
How to do this is given by the following (more difficult) exercise:

\begin{Exercise}[Rotation of coordinate system*]\label{Ex:rotation}
For given unit-vector $\vec u$ and angle $\theta$ define
$U=\exp(-i \theta \vec u\cdot\vec\sigma /2)$. Then $UU^*=U^*U=\mathbf 1$, i.e.,
$U$ is a unitary transformation of $\mathbb C^2$, and
$U\rho(\vec a)U^*=\rho(\vec b)$ where $\vec b \in \mathbb R^3$ results from
$\vec a$ by rotation about $\vec u$ through an angle $\theta$.
\end{Exercise}

This result really belongs to the representation theory of groups; 
a major topic having deep connections with quantum theory.
It is a curious fact that if $\theta=2\pi$ the operator $U$ is equal to
$-\mathbf 1$. So though $U$ works on a density matrix by a rotation through
$360^\circ$, it does not transform a state vector to itself but to its negative.
A rotation through $720^\circ$ or the angle $4\pi$ is needed to do this.
The fact that two complete revolutions are needed to transform a state
vector into itself whereas one revolution multiplies the state vector by
$-1$ has been experimentally verified through observation of interference
effects.

\subsection{Quantum Cram\'er-Rao inequality}
\label{ss1.3:quantum cramer-rao}
Consider a quantum statistical model whereby the density matrix $\rho$
depends on an unknown parameter $\theta$. Possibly $\theta$ is a vector
but we will not emphasize that fact in the notation. In particular,
a spin-half system has a density matrix $\rho=\rho(\vec a)$ depending on 
the vector $\vec a$ in the closed unit ball, which we will denote by $B$.
Interesting statistical models
could therefore have a one-, two- or three-dimensional parameter $\theta$,
specifying a curve, a surface, or an open region of $B$.
Of particular interest
are one- and two-dimensional {\em pure-state models\/} models, 
specifying a curve on the boundary $S$ of the unit sphere $B$ and 
the whole of $S$ respectively.
Results are strikingly different according to whether
the true value of $\theta$ corresponds to a point in $S$ or in the
interior of $B$. By a {\em mixed-state model\/} we mean a model in the
interior of $B$. By the {\em full model\/}, pure or mixed, we mean the model: `$\rho$ is in $S$', and `$\rho$ is in the interior of $B$' respectively.
By the natural parametrization of these models we mean the parametrization
$\rho=\rho(\vec u)$, $\rho=\rho(\vec a)$ respectively.\footnote{It 
would be nice to express conditions and results in the language of differential
geometry, i.e., independent of the specific parametrizations of the
models under consideration.}

The quantum Cram\'er-Rao bound involves a collection of self-adjoint matrices 
$\lambda_i$ called
the quantum score matrices, one for each component of $\theta$,
and a quantum information matrix. These
are defined as follows. 

\begin{Definition}[Quantum score matrices]
Suppose $\rho=\rho(\theta)$ depends on parameters
$\theta=(\theta_1,\dots,\theta_k)$. Suppose that $\rho$ is differentiable
with respect to $\theta$ and define self-adjoint matrices 
$\lambda_i=\lambda_i(\theta)$ implicitly by the equation
\begin{equation}
\label{e:quantum score}
\rho_i~=~\frac {\partial \rho(\theta)}{\partial\theta_i}~=~
\frac12(\lambda_i \rho + \rho \lambda_i).
\end{equation}
\end{Definition}

Note that the $\lambda_i=\lambda_i(\theta)$ will also depend on $\theta$.
Another name for these matrices is the symmetric logarithmic derivatives
of $\rho$ with respect to $\theta$. If $\rho$ and its derivative $\rho_i$ 
with respect to $\theta_i$ commute, then $\lambda_i$ is nothing else than
the derivative of $\log\rho$. By using a basis of $\mathbb C^d$ making $\rho$
diagonal, $\rho=\sum p_j \ket j \bra j$, one can solve \eqref{e:quantum score}
to obtain 
\begin{equation}\label{e:solve score}
\braopket j {\lambda_i} {j'} = \frac {2 \braopket j {\rho_i} {j'}}
{p_j+p_{j'}} .
\end{equation}
If some $p_j$ are zero the corresponding elements of $\lambda_i$ may be
chosen arbitrarily (subject to self-adjointness) without effect on subsequent
calculations. If $\rho$ is a pure state, then $\rho^2=\rho$ and it follows from
differentiating this equation with respect to $\theta_i$ that 
in this case $\lambda_i=2\rho_i$.

\begin{Exercise}[mean quantum score zero]\label{Ex:mean score zero}
Show that the quantum score has expectation zero, that is,
the distribution of a measurement of the observable $\lambda_i$ has
mean zero, or $\trace(\rho \lambda_i)=0$.
\end{Exercise}

\begin{Exercise}[Spin half, mixed]\label{Ex:spin half mixed}
Consider the full mixed-state spin-half model 
$d=2$, $\rho=\frac12(\mathbf 1+\vec\theta\cdot\vec\sigma)$, 
where $\theta$ is three-dimensional and satisfies $\sum_i\theta_i^2 < 1$. 
Then $\rho_i=\sigma_i$ for each $i$. At the point
$\theta=(0,0,\xi)$ the density matrix is diagonal with diagonal elements
$\frac12(1\pm\xi)$ and the quantum scores are found from \eqref{e:solve score}
to be
$\sigma_x$, $\sigma_y$ and $(1-\xi)^{-2}(-\xi{\mathbf 1}+\sigma_z)$.
\end{Exercise}

\begin{Exercise}[Spin half, pure]\label{Ex:spin half pure}
The full pure-state spin-half model has everything as in 
the previous exercise but now with $\sum_i\theta_i^2 = 1$.
A two-dimensional parametrization is called for, using, e.g., the polar coordinates of the unit
vector $\vec\theta$. However on the Northern hemisphere we can stick to
$\theta=(\theta_1,\theta_2)$ with $\theta_3=+(1-\theta_1^2-\theta_2^2)^{1/2}$
and we find that at $\theta=(0,0)$ the quantum scores are $\sigma_x$ and
$\sigma_y$.
\end{Exercise}

\begin{Exercise}[$n$ copies]
Suppose $\rho^{(n)}(\theta)=\rho^{\otimes n}(\theta)$.
Then the quantum scores are given by 
\begin{equation}\label{e:product score}
\lambda^{(1)}_i\otimes{\mathbf 1}\cdots\otimes{\mathbf 1}+\dots+
    {\mathbf 1}\otimes\cdots\otimes{\mathbf 1}\otimes\lambda^{(1)}_i.
\end{equation}
\end{Exercise}

Now we can define the quantum information matrix and state the original quantum 
Cram\'er-Rao bound.

\begin{Definition}[Quantum information matrix]\label{d:quantum information}
The quantum in\-for\-m\-a\-t\-ion matrix
$I_Q$ is defined by
\begin{equation}
(I_Q)_{ii'}~=~\frac12\trace(\rho(\lambda_i\lambda_{i'}+\lambda_{i'}\lambda_i)).
\end{equation}
\end{Definition}
Check that this defines a real, positive semi-definite matrix!

\begin{Exercise}[$n$ copies, continued]
Show from \eqref{e:product score} and Exercise \ref{Ex:mean score zero}
that the quantum information $I_Q^{(n)}$
for a parameter $\theta$ in the system $\rho^{\otimes n}$ is just $n$ times
the quantum information for $\theta$ in a single copy of the system.
\end{Exercise}

\begin{Theorem}[Quantum Cram\'er-Rao bound]
Define $I_M(\theta)$ to be the Fisher information matrix for
the parameter $\theta$ in the distribution of the outcome of a measurement
$M$ on the quantum system $\rho(\theta)$. Then (with respect to the
usual ordering of symmetric positive semi-definite matrices)
$I_M(\theta)\le I_Q(\theta)$.
\end{Theorem}
The result in this form was proved by Braunstein and Caves (1994) for
a one-dimensional parameter, but the general result is an easy consequence
of this by considering the information for arbitrary linear combinations.
As a corollary one obtains Helstrom's original form of the theorem as a
lower bound to the variance of an unbiased estimator of $\theta$ based on
the outcome of an arbitrary measurement $M$.

The proof is just as for the ordinary Cram\'er-Rao bound, an
exercise in using the Cauchy-Schwartz inequality, but now with the complex
inner-product $\trace X^* Y$ between two self-adjoint matrices. And just
as in the usual proof of the Cram\'er-Rao inequality, as a by-product
the proof shows that equality holds, for a one-dimensional
parameter, if (though not quite if and only if) $M$ is actually a simple
measurement of the observable $\lambda$:
\begin{Exercise}[Optimal $M$ for 1-d $\theta$]\label{Ex:optimal M 1d}
Show for one-dimensional $\theta$ that if $M$ is the 
simple measurement of the observable $\lambda$, i.e., its values are in
one-to-one correspondence with the eigenspaces of $\lambda$ and each
$M(x)$ is the projection onto the corresponding eigenspace, then
$I_M=I_Q$.
\end{Exercise}

There is a complication when using this result. Typically $\lambda$ will depend
on $\theta$, and typically in such a strong way that the eigenspaces of
$\lambda$ (and not just eigenvalues) depend on $\theta$. Thus the best
measurement of $\theta$ in terms of Fisher information depends on the
true value of $\theta$.
However things are very simple in the following example:
\begin{Exercise}Suppose all $\rho(\theta)$ commute, i.e.,
have common eigenspaces. Show that the
$\rho_i(\theta)$ then also commute for all $i$ and $\theta$. Show
that a simple measurement of the common eigenspaces of all these
matrices has Fisher information equal to the quantum information
for all values of $\theta$.
\end{Exercise}
The above is actually a completely classical model where 
$\rho=\sum_i p_i(\theta)\ket i \bra i$, i.e., a classical mixture 
with mixing distribution depending on $\theta$ of the fixed pure
states $\ket i$. The optimal measurement is to measure which of these
pure states the system is in; that can best be done using the 
projector-valued probability measure with elements $\ket i\bra i$
resulting in the outcome `$i$' with probability $p_i(\theta)$.
The quantum information matrix is the Fisher information matrix for
this distribution.

The result of Exercise \ref{Ex:optimal M 1d} does gives a lot of hope for 
a clear solution to the problem of
estimating a one-dimensional parameter, at least, for large $n$, for the 
system $\rho^{\otimes n}(\vec a(\theta))$, as was
first pointed out by Barndorff-Nielsen and Gill (1998). 
Suppose the parameter $\theta$ is identified, so
that there are a finite number of simple measurements, the distributions
of whose outcomes identifies $\theta$. For example, in the spin-half case
$\rho=\frac12(\mathbf 1 +\vec\theta\cdot\vec\sigma)$,
measurements of $\sigma_1$, $\sigma_2$, $\sigma_3$ result in Bernoulli
trials with probabilities $\frac 12(1\pm\theta_i)$. Suppose that from
consistent estimators of the $a_i(\theta)$ we can construct a consistent
estimator of $\theta$. Now, use a growing number but vanishing proportion 
of copies of our quantum system with which to `pre-estimate' $\theta$
consistently. Call this preliminary estimator $\widetilde \theta$. Now, compute the
quantum score for $\theta$ at $\widetilde\theta$, determine its
eigenspaces, and implement the corresponding simple measurement on all
remaining copies of the system. This gives us an i.i.d.\ sample from
some distribution $p(\cdot|\theta;\widetilde \theta)$. Estimate $\theta$
by maximum likelihood on these observations conditional on the observed
value of $\widetilde\theta$. The result $\widehat\theta$ will be an estimator
approximately normally distributed about $\theta$ with variance
approximately $1/nI(\theta;\widetilde\theta)$ where $I(\theta;\widetilde\theta)$
is the Fisher information for $\theta$ in one of these observations,
given $\widetilde\theta$. Now for $n$ large we have arranged that
$\widetilde\theta$ is close to the true value of $\theta$. We may
hope that the eigenspaces of $\lambda(\widetilde\theta)$ are close to
the eigenspaces of $\lambda(\theta)$ and hence that the Fisher information
in one observation (one simple measurement) of $\lambda(\widetilde\theta)$
is close to that in one observation of $\lambda(\theta)$. But the
latter achieves the quantum Cram\'er-Rao bound at $\theta$.
Thus under suitable smoothness conditions $I(\theta;\widetilde\theta)$ will 
be close to $I_Q(\theta)$ and hence the asymptotic distribution of our
final estimator close to normal about $\theta$ with variance 
$1/nI_Q(\theta)$. This is coming close to saying that $\widehat\theta$
is asymptotically optimal.

We know that no unbiased estimator of $\theta$ can 
have smaller variance. However that does not tell us no estimator whatever
can do better, e.g., in terms of mean square error. Indeed the phenomenon
of super-efficiency is just as present here as in ordinary statistics.
In order to make a compelling optimality statement about our estimator
we must either restrict attention to a sub-class of nicely behaved estimators,
or make optimality statements which are of a Bayesian or a
minimax nature. A very useful tool, which can be used in any of
these approaches, is the van Trees inequality which says for a one-dimensional
parameter $\theta$ with prior distribution $\pi(\d\theta)$, under some
regularity conditions, that the expected (with respect to the prior)
mean square error of a completely arbitrary estimator of $\theta$ is
bounded by one divided by the expected Fisher information for the parameter
plus the information, with respect to location, in the prior distribution.
This writer prefers to restrict the class of estimators according to some
regularity condition. We will go into this in more detail in the next section, but before that, let us
consider the multiparameter case. We will see that a more fundamental complication arises: at a fixed parameter value,
quantum scores for different components of the
parameter may not commute.

\begin{Exercise}[Quantum information for spin-half models]
In exercises \ref{Ex:spin half pure} and \ref{Ex:spin half mixed} we noted
the score matrices for the full pure-state model $\rho=\rho(\vec u)$
and for the full mixed-state model $\rho=\rho(\vec a)$.
Show that, at $\vec u=(0,0,1)$ in the first case and at $\vec a=(0,0,\xi)$
in the second case, the quantum information matrices for $\theta=(u_1,u_2)$
and for $\theta=\vec a$ are respectively
\begin{equation}
I_Q=\left(
\begin{array}{cc}
1 & 0 \\
0 & 1
\end{array}
\right),
\qquad
I_Q=\left(
\begin{array}{ccc}
1 & 0 & 0\\
0 & 1 & 0\\
0 & 0 & 1/(1-\xi)^2
\end{array}
\right).
\end{equation}
\end{Exercise}

Now the approach just sketched for the one-parameter case breaks down.
Certainly we can form a preliminary estimator of $\theta$ and thereby
`estimate' the quantum score matrices. Next, in the full pure- and mixed-state
models, one can rotate the coordinate system and reparametrize 
so that the quantum scores become
$\sigma_x,\sigma_y$ (pure-state model), and
$\sigma_x,\sigma_y, a+b\sigma_z$ (mixed state model).
There is no way, in either model,
we can simultaneously measure these observables since they
do not commute. Thus no measurement on a single particle has an information
matrix equal to $I_Q$. The big question is, what {\em is\/} the class of
information matrices $I_M$ which are available? And if we
can perform measurements on the system obtained by combining particles,
what scaled information matrices $I_M^{(n)}/n$ become available? The
latter class includes all of the former class,
since the joint measurements include $n$ i.i.d.\ copies of measurements on
separate particles; moreover these classes are convex and bounded.

Though all scaled information matrices $I_M^{(n)}/n$ are bounded by
$I_Q$, we cannot expect them, for given $n$, to contain a single
`best' information. Which measurement we should choose will depend on
the relative accuracy with which we want to estimate the different components
of $\theta$. For instance if in the pure-state case, close to 
$\theta=(0,0)$, we are only interested in $\theta_1$ we should simply
measure $\sigma_x$ on each of our $n$ particles yielding the maximum
information on $\theta_1$ but no information at all on $\theta_2$.
After we have characterized the class of all information matrices available,
we must specify through some loss function the relative importance of
the different parameters and solve some optimization problem.

\section{A new Cram\'er-Rao type bound}
In this section we report on recent results of Gill and Massar
(1998), concentrating on the spin-half situation, and within that case,
emphasizing the full pure-state model and the full mixed-state model.
There turns out to be a striking difference between these two cases.
For pure states, there is asymptotically no advantage
in joint measurements on many particles. However for mixed states there
typically is an advantage. How much is still an open question.
The following result should be called a `Theorem' (in quotes) 
since we do not specify regularity conditions and indeed only a `Proof'
exists, not yet a Proof.

\begin{Theorem}[Achievable information matrices, $n=1$]\label{first main theorem}
The set of all information matrices of outcomes of measurements
of one spin-half particle for a smooth model $\rho(\theta)$
is $\{F:\trace(I_Q^{-1}F)\le 1\}$.
\end{Theorem}

The parameter $\theta$ could be one-, two- or three-dimensional.
We suppose either that we have a pure-state model, or a strictly 
mixed-state model.  The argument, in Gill and Massar (1998), has
two main parts. In the first part we show that
for all $M$, $F=I_M$ satisfies $\trace(I_Q^{-1}F)\le d-1$
(we do not yet need that $d=2$).
In the second part we show that, when $d=2$, for
any $F$ satisfying this inequality one can construct a measurement $M$
for which $I_M=F$. For $d>2$, not all
$F$ satisfying $\trace(I_Q^{-1}F)\le 1$ are achievable, and it remains open
to characterize exactly the class of achievable information matrices.

For the first part a series of preparatory steps are taken to bring us,
`without loss of generality', to a
situation that allows exact computations. For simplicity take $d=2$.
If $\rho(\theta)$ lies in the interior of the unit ball,
and $\theta$ has dimension one or two, one can augment $\theta$ with other
parameters, raising its dimension to $3$. This can be done in such a way
that the cross-information elements in the augmented $I_Q(\theta)$ are
all zero. It then suffices to prove the inequality for $\theta$ of dimension
$3$, and then we may as well use the natural parametrization $\rho(\theta)=
\frac12(\mathbf 1 +\vec\theta\cdot\vec\sigma)$ 
with $\|\theta\|<1$  since the the quantity 
$\trace(I_Q^{-1}F)$ is invariant under smooth reparametrization.
If on the other hand $\rho(\theta)$ is a pure state model we can in the
same way after augmenting $\theta$ assume that $\theta$ has dimension $2$
and after reparametrization the model is  $\rho(\theta)=
\frac12(\mathbf 1 +\vec\theta\cdot\vec\sigma)$ with $\|\theta\|=1$.

For the next preparatory step we need the concepts of refinement and 
coarsening of a measurement. 
\begin{Definition}[Coarsening and refinement] A measurement
$M$ with sample space $\mathcal X$ is a refinement of 
$M'$ with sample space $\mathcal Y$ (and $M'$ is a coarsening of $M$)
if a measurable function $f:\mathcal X \to \mathcal Y$ exists with
$M'(B)=M(f^{-1}(B))$.
\end{Definition}
The result of measurement of $M'$ then has the same
distribution as taking $f$ of the outcome of measurement of $M$.
It follows that the Fisher information in the outcome of $M'$ is
less than or equal to that in $M$ since under coarsening of data,
Fisher information can only decrease. 

Now we show that any measurement $M'$ has a refinement $M$ for
which $M(A)=\int_A M(x)\mu(\d x)$ for some nonnegative operator-valued
function $M$ and bounded measure $\mu$ 
{\em and for which $M(x)$ has rank one for all $x$\/}, thus $M(x)=\ket{\psi(x)}\bra{\psi(x)}$ for some 
(not necessarily normalised) vector function $\ket{\psi(x)}$.
Consequently it will suffice to prove the result for such 
{\em maximally refined\/} measurements $M$.
Start with the measurement $M'$ with sample space $\mathcal Y$.
Define a probability measure $\nu$ on $\mathcal Y$ by
$\nu(B)=\trace(M'(B))/d$; 
by taking Radon-Nikodym derivatives one can define $M'(y)$
such that $M'(B)=\int_B M'(y)\nu(\d y)$. Since the
rank of $M'(y)$ is finite, $M'(y)=\sum_i M_i(y)$ where
each $M_i(y)$ has rank one. Now refine the original sample space $\mathcal Y$
to $\mathcal X =\mathcal Y \times \{1,\dots,d\}$, defining $M(A\times\{i\})=
\int_A M_i(y)\nu(\d y)$. Equivalently $M(A)=\int_A M_i(x) \mu(\d x,\d i)$
where $\mu$ is the product of $\nu$ with counting measure.

This brings us to the situation where the model is either full pure-state 
or full mixed-state, and where the measurement is maximally
refined. We take the natural parametrization of either of these models,
and without loss of generality work at a point $\theta$ where $\theta=(0,0)$
or $(0,0,\xi)$. This is possible by the result of Exercise \ref{Ex:rotation}.
Now we have a formula for $I_Q$ and for the derivatives of $\rho$
with respect to the components of $\theta$, both in the pure and the mixed
case, and we have a representation
for $M$ in terms of a collection of vectors $\psi(x)$ which must satisfy
the normalization constraint 
$\int_{\mathcal X} \ket{\psi(x)}\bra{\psi(x)}\mu(\d x)=\mathbf 1$ but which are otherwise
arbitrary. Both $\rho$ and $I_Q$ are diagonal. We simply compute
$\trace I_Q^{-1} I_M$ and show that it equals $1$ in the case $d=2$.
We leave the details as an exercise for the diligent reader---the 
computation is not difficult but does not seem all that illuminating either. 
We would dearly like to know if there a more insightful way to get this result!

The same arguments work for arbitrary $d$ though the details are more
complicated; a full mixed-state model has $\frac12 d (d+1)$
parameters, a full pure-state model $\frac12 d (d+1) -(d-1)$
parameters, and a careful parametrization is needed to make $I_Q$
diagonal.

In the second part (for $d=2$ only) it is shown that for
any $F$ satisfying $\trace(I_Q^{-1}F)\le 1$, one can construct a 
measurement $M$ for which $I_M=F$. This measurement will be described
in the next section. It typically depends on the point $\theta$ so
a multi-stage procedure is going to be necessary to achieve 
asymptotically this information bound. That will be the main content
of the next section, where we do some quantum asymptotics proving 
asymptotic optimality results for $n\to\infty$ of the resulting
two-stage procedure.

We only have partial results for $n>1$. In two special cases the
available scaled information matrices do not increase as $n$ increases.
One of these cases is the case of pure-state models. This case 
has been much studied in the
literature and is of great practical importance. The other case is when
we make a restriction on the class of measurements to measurements of
product form (in the literature also sometimes called an
{\it unentangled measurement\/}). We first define this notion and then explain its significance.

\begin{Definition}[Product-form measurements]
We say that a measurement on $n$ copies of a given quantum system is
of product form if $M^{(n)}(A)=\int_A M^{(n)}(x)\mu(\d x)$ for a 
real measure $\mu$
and matrix-valued function $M^{(n)}(x)$ where 
$M^{(n)}(x)$ is of the form $M_1(x)\otimes\cdots\otimes M_n(x)$, with
nonnegative components.
\end{Definition}

We described in the previous section a measurement procedure whereby we
first carried out measurements on some of our $n$ particles, and then
depending on the outcome, carried out other measurements on the remaining
particles. Altogether this procedure constitutes one measurement on the
joint system of $n$ particles taking values in some $n$-fold product
space. One can conceive of more elaborate schemes where depending on
the results at any stage, one decides, possibly with the help of some
classical randomisation, which particle to measure next and how. It
would be allowed to measure again a particle which had previously been
subject to measurement. There exists a general description of the state
of a quantum system after measurement, allowing one to piece all the
ingredients together into one measurement of the combined system.
A measurement which can be decomposed into separate steps
consisting of measurements on separate particles only, is called a
{\em separable\/} measurement.

It turns out that all separable measurements (provided all
outcomes of the component steps are encoded in the overall outcome
$x$) have product-form.
On the other hand, product-form measurements exist which are
not separable, see Bennett et al.\ (1998). The product-form measurements
form a large and interesting class, including all measurements which can 
be carried out sequentially on separate particles as well as more besides.

In the notion of separable measurement it is insisted that all
intermediate outcomes are included in the final outcome. If one throws
away some of the data, one gets an outcome whose distribution is the
same as the distribution of a coarsening of the original measurement.
Coarsening of a measurement can easily destroy the properties of being
separable or being of product-form. This is some explanation for
the complicated restriction to measurements which can be refined to
product-form in the following theorem:

\begin{Theorem}[Achievable information matrices $n>1$]\label{second main theorem}
The scaled information
matrices of measurements on a smooth model $\rho^{\otimes n}(\theta)$
remain $\{F:\trace(I_Q^{-1}F)\le 1\}$
\begin{enumerate}
\item in a pure-state spin-half model;
\item in a mixed-state spin-half model with the class of measurements
 restricted to measurements which can be refined to product-form, 
\end{enumerate}
\end{Theorem}
The theorem is proved exactly as before, again finishing in an
unilluminating calculation.

We have a counterexample to the conjecture that, for mixed states,
the bound holds
for {\em all\/} measurements. In
the case $n=2$, at the point $\rho=\frac12\mathbf 1$, there is
a measurement for which $\trace(I_Q^{-1}I_M^{(2)}/2)=3/2$,
thus $50\%$ more information in an appropriate measurement of two
identical particles than any combination of separate measurements of
the two. What the set of achievable scaled information matrices looks
like and whether it continues to grow (and to what limit) as $n$ grows
is completely unknown.

The measurement has seven elements, the first six
of the form $\frac12\Pi_{[\psi]}$, and the seventh $\Pi_{[\psi]}$. 
The various $\psi$ are 
$\ket{+z+z}$, $\ket{-z-z}$, $\ket{+x+x}$, $\ket{-x-x}$, 
$\ket{+y+y}$, $\ket{-y-y}$, $\ket{S}$. By $\ket{+z+z}$ we
mean $\ket{+z}\otimes\ket{+z}=\psi(\vec e_z)\otimes\psi(\vec e_z)$ and similarly for the next five.
The last $\psi$ is the so-called {\em singlet state\/}
$\frac 1{\sqrt{2}}(\ket{+z}\otimes\ket{-z}-\ket{-z}\otimes\ket{+z})$.
As a pure state of two interacting spin-half particles, this is the
famous entangled state resulting in the violation of the Bell inequalities,
and hence of locality (according to some interpretations).
Here it arises as part of a measurement of two completely non-interacting
particles; however this measurement can never be implemented by doing
separate operations on the separate particles.

Similar examples occur in the paper of Vidal et al.\ (1998), extending
the pure-state results of Massar and Popescu (1995) to mixed states.

\section{Quantum asymptotics}
The results of the previous section are in the form of a bound
on the information matrix based on the outcome of any measurement (perhaps
restricted to the class of product-form measurements) 
on $n$ identical copies of a given spin-half
quantum system with state depending on an unknown parameter $\theta$.
We will now explain how such a bound can be used to give
asymptotic bounds on the quality of estimators based on those measurements.
Furthermore, we
show how the bounds can be achieved by a two-stage procedure
using simple measurements on separate particles only.
As far as achieving the bounds is concerned, only for the full
mixed-state model under the natural parametrization 
is the problem completely solved. For the other models, the results are
conjectural. 

We will discuss two kinds of bounds: firstly, a bound on the limiting 
scaled mean quadratic error matrix of a well-behaved sequence of estimators,
and secondly, a bound on the mean quadratic error matrix of the limiting
distribution of a well-behaved sequence of estimators. Each
has its advantages and disadvantages. In particular, since the delta-method
works for (the variance of) limiting distributions but not for 
limiting mean square errors, stronger conditions are needed to prove optimality
of some procedure in the first sense than in the second sense.

\subsection{Two asymptotic bounds}
Obviously a bound on the information matrix, by the ordinary Cram\'er-Rao
inequality, immediately implies a bound on the covariance matrix of an
unbiased estimator. However this is not a restriction we want to make.
It turns out much more convenient to work via a Bayesian version of the
Cram\'er-Rao inequality due to van Trees (1968), as generalised to the
multi-parameter case by Gill and Levit (1995). For a one-dimensional
parameter the van Trees inequality is easy to state: the Bayes quadratic
risk is bounded by one over expected information plus information in the prior.
In the multiparameter case one has a whole collection of inequalities 
corresponding to different choices of quadratic loss function and 
some other parameters, more difficult to interpret.

Let $\pi(\theta)$ be a prior density for the $p$-dimensional parameter $\theta$,
which we suppose to be sufficiently smooth and supported by a compact 
and smoothly bounded region of
the parameter space; see Gill and Levit (1995) for
the precise requirements. Let $C(\theta)$ be a $p\times p$ symmetric
positive definite matrix (C stands for cost function) 
and let $V_M^{(n)}(\theta)$ be the mean quadratic
error matrix of a chosen estimator of $\theta$ based on a measurement of $n$
copies of the quantum system. Letting $\Theta$ denote a random drawing from
the prior distribution $\pi$, it follows that $\E\trace C(\Theta)V_M^{(n)}(\Theta)$
is the Bayes risk of the estimator with respect to the loss function
$(\widehat\theta^{(n)} - \theta)^\top C(\theta) 
(\widehat\theta^{(n)} - \theta)$.

Let $D(\theta)$ be another $p\times p$ matrix function of $\theta$.
Let $I_M^{(n)}(\theta)$ denote the Fisher information matrix in
the measurement. Then the multivariate van Trees inequality reads
\begin{equation}\label{van trees}
\E\trace C(\Theta)nV_M^{(n)}(\Theta) \ge
\frac{ (\E\trace D(\Theta))^2 }
{\E\trace C(\Theta)^{-1}D(\Theta)(I_M^{(n)}(\Theta)/n)D(\Theta)^\top +
 \widetilde {\mathcal I}(\pi)/n }
\end{equation}
where
\begin{equation}\label{I tilde}
\widetilde {\mathcal I}(\pi) = \int\frac 1 {\pi(\theta)}
\sum_{ijkl}C_{ij}^{-1}(\theta)
\frac \partial {\partial\theta_k}
( D_{ik}(\theta)\pi(\theta) )
\frac \partial {\partial\theta_l}
(D_{jl}(\theta)\pi(\theta) )
\d\theta.
\end{equation}
On invoking Theorem \ref{second main theorem} we have the bound
$\trace I_Q^{-1}(\theta)(I_M^{(n)}(\theta)/n)\le 1$,
 (provided that, in the mixed case, we
restrict attention to measurements refinable to product-form). 
We are going to assume
that our sequence of measurements and estimators is such that the
normalized mean quadratic error matrix $V_M^{(n)}(\theta)$ converges
sufficiently regularly to a limit $V(\theta)$. Our aim is to transfer
the just mentioned bound to $V$ obtaining the bound 
$\trace I_Q^{-1}(\theta)V(\theta)^{-1}\le 1$.

We will do this by making appropriate choices of $C$ and $D$. 
We will need regularity conditions both on the sequence of
estimators and on the model $\rho(\theta)$ in order to carry
over equation \eqref{van trees} to the limit.

\begin{Theorem}[Asymptotic Cram\'er-Rao 1]\label{third main theorem}
Suppose that on some open set of parameter values $\theta$:
\begin{enumerate}
\item $nV_{(n)}$ converges uniformly
to a continuous limit $V$. 
\item $I_Q(\theta)$ is continuous with bounded partial derivatives.
\item $V$ and $I_Q$ are non-singular.
\end{enumerate}
Then the limiting normalised mean quadratic error matrix satisfies
\begin{equation}\label{asymptotic c-r}
\trace I_Q^{-1}(\theta)V(\theta)^{-1}\le 1.
\end{equation}
\end{Theorem}

We outline the proof of the theorem as follows. First of all,
we pick a point $\theta_0$ and define $V_0=V(\theta_0)$. Next
we define
\begin{equation}
C(\theta)=V_0^{-1}I_Q^{-1}(\theta)V_0^{-1},
\end{equation}
\begin{equation}
D(\theta)=V_0^{-1}I_Q^{-1}(\theta).
\end{equation}
With these choices \eqref{van trees} becomes
\begin{equation}\label{van trees 2}
\E\trace V_0^{-1}I_Q^{-1}(\Theta) V_0^{-1} (nV_M^{(n)}(\Theta)) \ge
\frac{ (\E\trace V_0^{-1} I_Q^{-1}(\Theta))^2 }
{\E\trace I_Q(\Theta)^{-1}(I_M^{(n)}(\Theta)/n) +
 \widetilde {\mathcal I}(\pi)/n }.
\end{equation}
We can bound the first term in the denominator of the right hand side
by $1$, by the results of the last section. The second term in the
denominator of the right hand side is finite, by our third assumption,
and for $n\to\infty$ it converges to zero. By our first 
assumption \eqref{van trees 2} converges to
\begin{equation}\label{van trees 3}
\E\trace V_0^{-1}I_Q^{-1}(\Theta) V_0^{-1} V(\Theta) \ge
 (\E\trace V_0^{-1} I_Q^{-1}(\Theta))^2 .
\end{equation}
Now replace the prior density $\pi$ by one in a sequence of priors,
concentrating on smaller and smaller neighbourhoods of $\theta_0$.
Using the continuity assumptions on $V$ and $I_Q$, we obtain from
\eqref{van trees 3} the inequality
\begin{equation*}
\trace V_0^{-1}I_Q^{-1}(\theta_0) V_0^{-1} V_0 \ge
 (\trace V_0^{-1} I_Q^{-1}(\theta))^2 .
\end{equation*}
or in other words, with $\theta=\theta_0$, the required
\begin{equation}
\trace I_Q^{-1}(\theta)V^{-1}(\theta) \le
 1.
\end{equation}

In some situations it might be more convenient to have a bound on the 
mean quadratic error of a limiting distribution, assuming one to exist.
At the moment of writing we believe the following:

\begin{Theorem}[Asymptotic Cram\'er-Rao 2]\label{fourth main theorem}
Suppose
\begin{enumerate}
\item $\widehat \theta_n$ is H\'ajek regular at $\theta$ at root $n$ rate.
\item If $Z$ has the limiting distribution of $\sqrt{}n(\widehat\theta-\theta)$,
then the mean quadratic error matrix of the limiting distribution 
$V=\E(Z\,Z^\top)$ is non-singular.
\item $I_Q$ is non-singular.
\end{enumerate}
Then $V$ satisfies
\begin{equation}\label{asymptotic c-r 2}
\trace I_Q^{-1}(\theta)V(\theta)^{-1}\le 1.
\end{equation}
\end{Theorem}
The proof should follow the lines of the similar result in Gill and Levit (1995), with a prior distribution concentrating on a root $n$ neighbourhood
of the truth. We will need similar choices of $C$ and $D$ as in the proof of
Theorem \ref{third main theorem} though the dependence of $D$ on $\theta$
can now be suppressed.

\subsection{Achieving the asymptotic bounds}
At present we have essentially complete results in the full mixed-state
spin-half model with the natural parametrization. We believe they
can be extended to smooth ($C^1$) pure- and mixed-state models.

Give yourself a target mean quadratic error matrix $W(\theta)$
satisfying 
\begin{equation}\label{constraint}
\trace I_Q(\theta)^{-1}W(\theta)^{-1}\le 1.
\end{equation}
Is there a sequence of measurements $M^{(n)}$ satisfying the
conditions of Theorems \ref{third main theorem} or \ref{fourth main theorem}
with limiting
mean quadratic error matrix $V(\theta)$ equal to the target?

Possibly we do not start with a target $W$ but with a
step earlier, with a quadratic cost function.
For given $C(\theta)$ it is straightforward to compute the matrix
$W(\theta)$ which minimizes $\trace C(\theta)W(\theta)$
subject to the constraint \eqref{constraint}; the solution is
$W=\trace(( I_Q^{-\frac12} C I_Q^{-\frac12})^{\frac12})
I_Q^{-\frac 12} ( I_Q^{\frac12} C I_Q^{\frac12})^{\frac12} I_Q^{-\frac12}$.
Now we pose the same question again, with the $W$ we have just calculated
as target. 

Let us call $F=W^{-1}$ the target information matrix.
First we pretend $\theta$ is known and exhibit a measurement $M$ on
a single particle with the target information matrix at the given
parameter value.

In the previous section we omitted explaining how the bound of Theorem
\ref{first main theorem} can be attained. That theorem stated that,
at a given parameter value, for any positive-semidefinite symmetric $F$ 
satisfying 
$\trace I_Q^{-1}F\le 1$ there is a measurement $M$ on a single spin-half 
particle with $I_M=F$.  What is that measurement? We describe it in the
case of a full mixed-state spin-half model with the natural parametrization, 
thus $\rho(\theta)=\frac 12(\mathbf 1
+\vec\theta\cdot\vec\sigma)$. The matrices $I_Q$ and $F$ are $3\times 3$.

To start with, we compute the eigenvector-eigenvalue decomposition
of $I_Q^{-\frac12} F I_Q^{-\frac 12}$, obtaining eigenvectors $\vec h_i$
and nonnegative eigenvalues $\gamma_i$, say. The condition on $F$ translates to
$\sum\gamma_i\le 1$. Now define $\vec g_i=I_Q^{\frac 12}\vec h_i$
and three unit vectors
$\vec u_i= g_i   /\| g_i/ \|$,
and finally consider the measurement $M$ taking seven different values, whose
elements are $\gamma_i\Pi(\pm\vec u_i)$,  $i=1,2,3$, and
$(1-\sum\gamma_i)\mathbf 1$.

It turns out by a staightforward computation (carried out, without
loss of generality, at $\theta=(0,0,\xi)$) that the information matrix
for the measurement with the two elements $\Pi(\pm \vec u_i)$ has
information matrix $\vec g_i\otimes \vec g_i$ 
and hence the measurement $M$
has information matrix $\sum_i \gamma_i\vec g_i\otimes\vec g_i=F$.

This seven-outcome measurement can be implemented as a randomized choice 
between three simple measurements: with probability
$\gamma_i$ measure spin in the direction $\vec u_i$, with probability
$1-\sum \gamma_i$ do nothing.

However in practice this measurement is not available since the directions
$\vec u_i$ and probabilities $\gamma_i$ depend on the unknown $\theta$.
We therefore take recourse to the following two-stage measurement procedure.

First measure spin in the $x$, $y$ and $z$ directions on $\frac 13 n^a$
each of the particles, where $0<a<1$ is fixed and the numbers are rounded
to whole numbers. The expected relative frequency of `up' particles in each 
direction
is $\frac 12 (1+\theta_i)$, $i=1,2,3$, so solving observed equals expected
yields a consistent preliminary estimator $\widetilde\theta$ of $\theta$.
If the estimate lies outside the unit-ball project onto the ball and
stop. With large probability no projection is necessary.
We can compute the eigenvalue-eigenvector decomposition of
$I_Q^{-\frac12}(\widetilde\theta) 
F(\widetilde\theta) I_Q^{-\frac 12}(\widetilde\theta)$, leading to
fractions $\gamma_i$ and directions $\vec u_i$ as above. Measure the spin
of a fraction $\gamma_i$ of the remaining particles in the direction
$\vec u_i$. Solve again the three (linear) equations 
`observed relative frequency
equals expected' treating the $\vec u_i$ as fixed. Project onto
the unit ball if necessary, yielding an estimator $\widehat \theta$.

Our claim is that this procedure exhibits a measurement $M^{(n)}$
on the $n$ particles, and an estimator $\widehat \theta^{(n)}$
based on its outcome, which satisfies the conditions of 
Theorem \ref{third main theorem}, with $V(\theta)$ equal to the
target $W(\theta)$. Thus the bound of Theorem \ref{third main theorem}
is also achievable, and a measurement which does this has been 
explicitly described above. Apart from projecting onto the unit ball
the estimator involves only linear operations on binomial variables
so is not difficult to analyse explicitly.
We need a preliminary sample size $\widetilde n$ of order $n^a$ and not, 
for example, of order $\log n$, in order to control the scaled mean quadratic error of the 
estimator. There is an exponentially small probability---in $\widetilde n$, not in $n$---that the preliminary estimate is outside of a given neighbourhood
of the truth, and hence that the scaled quadratic error is of order $n$.

One can further check that the estimator we have described also satisfies
the conditions of Theorem \ref{fourth main theorem}. 

Possibly one is interested in a different parametrization of the model.
Under a smooth ($C^1$) reparametrization, the delta method allows us to maintain
optimality in the sense of Theorem \ref{fourth main theorem}. However 
optimality in the sense of Theorem \ref{third main theorem} could be destroyed;
in order for it to be maintained the reparametrization should also be
bounded. Alternatively one must
modify the estimator by a truncation at a level increasing slowly 
enough to infinity with $n$, cf.\ Schipper (1997; section
4.4) or Levit and Oudshoorn (1993) for examples of the technique.

This approach can be extended to other spin-half models.
The difficulties are exemplified by the case of the two-parameter
full pure-state spin-half model.
Locally, consider the natural parametrization $\theta=(\theta_1,\theta_2)$,
$\theta_3=(1-\theta_1^2-\theta_2^2)^{1/2}$, $\rho=\rho(\vec\theta)$
at the point $\theta=(0,0)$. The quantum information
matrix for three parameters $\theta_1,\theta_2,\theta_3$
contains an infinite element.
However the recipe outlined above continues to work if we add to
a given $2\times 2$ target information matrix a third zero row and 
column---infinities always get multiplied by zero. The third fraction
$\gamma_3=0$ so simple measurements in just two directions suffice.

The resulting procedure involves linear operations on binomial counts,
projecting onto $S$, and reparametrization. Under some smoothness we
should finish with an estimator optimal in the sense of Theorem
\ref{fourth main theorem}; under further smoothness, boundedness,
and a sufficiently large preliminary sample also optimality in the
sense of Theorem \ref{third main theorem} should hold.

If the target information matrix includes some zeros, i.e., one
is not interested at all in certain parameters, the results should still
go through; the preliminary sample should be of size of order $n^a$, 
$\frac12 <a<1$, in order that the uncertainty in the initial estimate
of the `nuisance parameters' does not contaminate the final result.

\section{Non-locality without entanglement}
It would take us too far afield here to explain the notions of entanglement and of
non-locality. For some kind of introduction see K\"ummerer and Maassen (1998)
and Gill (1998), and Gill (1995a, 1995b); see also the books of Peres (1995),
Isham (1995), Penrose (1994), Maudlin (1994). However we would like to discuss whether or not
our finding, that non-separable joint measurements on
several independent (non-entangled) quantum particles can yield more
information that any separate measurements on the separate particles,
should be considered surprising or not. Recall that separable measurements,
cf.\ Bennett et al.\ (1998), 
are measurements which can be decomposed into a sequence of measurements
on separate particles, each measurement possibly depending on the outcome of
the preceding ones, and whereby it is allowed to measure further a
particle which has already been measured (and hence its state has been
altered in a particular way) at an earlier step.

>From a mathematical point of view there should not be much surprise.
The class of separable measurements is contained in the class of
product-form measurements, which is clearly a very small part of the
space of all measurements whatsoever. The optimisation problem
of maximising Fisher information (more precisely, some scalar functional 
thereof) must only be expected to have a larger outcome when we optimise
over a larger space. The surprise for the mathematician is rather that 
for pure states, and for one dimensional parameters, there is no gain in 
joint measurements. And it is strange that mixed states should exhibit
this phenomenon whereas pure states do not: the differenence is classical
probabilistic mixing which should not lead to nonclassical behaviour.

However physicists are and should be surprised.
The reason is connected to the feeling of many physicists
that the randomness in measurement
of a quantum system should have a deterministic explanation
(Einstein: ``God does not throw dice'') . We appreciate
very well that tossing a coin is essentially a completely deterministic
process. It is only uncontrolled variability in initial conditions which
lead to the outcome appearing to be completely random. Might it be the
case also that the randomness in the outcome of a measurement of a
quantum system might be `merely' the reflection of statistical variability
in some initial conditions? So-called hidden variables because at present
no physicist is aware what these lower level variables are and there is no
known way directly to measure them?

In fact there already exist arguments aplenty that if there {\it is\/} a
deterministic hidden layer beneath quantum theory, it violates
other cherished physical intuitions, in particular the principle of
locality; see again K\"ummerer and Maassen (1998), Gill (1998) for
some introduction to the phenomenon of entanglement, and further references.
But let us ignore that evidence and consider the new evidence from the
present results. Consider two identical copies of a given quantum state.
Suppose there were a hidden deterministic explanation for the randomness
in the outcome of any measurement on either or both of these particles.
Such an explanation would involve hidden variables $\omega_1$, $\omega_2$
specifying the hidden state of the two particles. Since applying
separate measurements to the two systems produces independent outcomes,
and since the outcomes of the same measurements are identically distributed,
one would naturally suppose that these two variables are independent and
identically distributed. Their distributions would of course depend
on the unknown parameter $\theta$. Now when we measure the joint system,
there could be other sources of randomness in our experiment, possibly
even quantum randomness, but still it would not have a distribution 
depending on $\theta$. So let us assume there is a third random
element $\omega_M$ such that the outcome of the measurement $M$ on
the system $\rho(\theta)\otimes\rho(\theta)$ is a deterministic function
of $\omega_1$, $\omega_2$ and $\omega_M$; the first two are independent
and identically distributed, with marginal distributions depending on
$\theta$, while the distribution of $\omega_M$ given the other two is
independent of $\theta$. Thus the random outcome $X$ of the measurement of $M$
is just $X(\omega_1,\omega_2,\omega_M)$, a random variable on the 
probability space
$(\Omega\times\Omega\times\Omega_M), ((\P_\theta\times\P_\theta)*\P_M)$ where
$\P_M$ is some Markov kernel from $\Omega\times\Omega$ to $\Omega_M$. Now
it is well-known from ordinary statistics that the Fisher information in
$\theta$ from the distribution of any random variable defined on this
space is less than twice the information in one observation of $\omega_1$
itself seen as a random variable defined on $(\Omega,\P_\theta)$.
Thus if one could realise any $\Omega_M,\P_M$ and any $X$ whatsover by
suitable choice of measurement $M$, achievable Fisher information would be additive!

What can we conclude from the fact that achievable Fisher information is not additive? We cannot rule out hidden variable models such as the above. But
apparently, the hidden variables are so well hidden that we cannot
uncover them from any measurements on single particles. i.e., it
is not possible to realise any $(\Omega_M,\P_M)$ and any $X$ whatever by
appropriate choice of experimental set-up. However we can uncover
the hidden variables better, apparently, from appropriate measurements 
on several particles brought together, even though these particles have
nothing whatever to do with
one another---their hidden variables are independent and identically
distributed. Alternatively the explanation must be found in 
some pathological non-measurability
or non-regularity of the statistical model we have just introduced.
Whatever escape-route one chooses, it is clear that if there is a 
deterministic explanation for quantum randomness, it is a very very weird explanation. God throws rather peculiar dice.

\acknowledgements
{This paper is based on work in progress together with
O.E. Barndorff-Nielsen and with S. Massar. I am grateful for the hospitality
of the Department of Mathematics and Statistics, University of Western
Australia. I would like to thank Boris Levit for his patient advice.}

\references

\Ref{
Barndorff-Nielsen, O.E. and Gill, R.D. (1998). 
{\it An example of non-attainab\-i\-l\-ity of expected quantum information}. 
\newline Preprint {\tt quant-ph/9808009}, {\tt http://xxx.lanl.gov}.
}

\Ref{
Bennett, C.H., DiVincenzo, D.P., Fuchs, C.A., Mor, T., Rains, E., Shor, P.W., Smolin, J.A., and Wootters, W.K. (1998).
{\it Quantum nonlocality without entanglement}.
\newline Preprint {\tt quant-ph/9804053}, {\tt http://xxx.lanl.gov}.
}

\Ref{
Biane, P. (1995). Calcul stochastique non-commutatief. pp.\ 4--96
in: {\it Lectures on Probability Theory: Ecole d\'et'e de Saint Flour
XXIII--1993}, P. Biane and R. Durrett, Springer Lecture Notes in
Mathematics {\bf 1608}.
}

\Ref{
Braunstein, S.L. and Caves, C.M. (1994).
Statistical distance and the geometry of quantum states.
{\it Physical Review Letters \bf 72}, 3439--3443.
}

\Ref{
Brody, D.C. and Hughston, L.P. (1998), Statistical geometry in quantum 
mechanics,
{\it Proceedings of the Royal Society of London Series A \bf 454}, 2445--2475.
}

\Ref{
Gill, R.D. (1995a). {\it Discrete Quantum Systems}.
Lecture notes \newline {\tt http://www.math.uu.nl/people/gill/discrete.ps.gz}.
}

\Ref{
Gill, R.D. (1995b). {\it Notes on Hidden Variables}.
Lecture notes \newline {\tt http://www.math.uu.nl/people/gill/hidden.ps.gz}.
}

\Ref{
Gill, R.D. (1998). Critique of `Elements of quantum probability'.
{\it Quantum Probability Communications \bf 10}, 351--361.
}

\Ref{
Gill, R.D. and Levit, B.Y. (1995).
Applications of the van Trees inequality: a Bayesian Cram\'er-Rao bound.
{\it Bernoulli \bf 1} 59--79
}

\Ref{
Gill, R.D. and Massar, S. (1998).
State estimation for large ensembles.
Preprint \newline {\tt http://www.math.uu.nl/people/gill/temp/massar7.ps}.
}

\Ref{
Helstrom, C.W. (1976). 
{\it Quantum Detection and Estimation Theory}. 
Academic, New York.
}

\Ref{
Holevo, A.S. (1982). 
{\it Probabilistic and Statistical Aspects of Quantum Theory}. 
North Holland, Amsterdam.
}

\Ref{ 
Holevo, A.S. (1983). Bounds for generalized uncertainty of the shift parameter. {\it Springer Lecture Notes in Mathematics \bf 1021}, 243--251.
}

\Ref{
Isham, C. (1995). {\it Quantum Theory}.
World Scientific, Singapore.
}

\Ref{
K\"ummerer, B. and Maassen, H. (1998). Elements of quantum probability.
{\it Quantum Probability Communications \bf 10}, 73--100.
}

\Ref{
Levit, B.Y.  and Oudshoorn, C.G.M. (1993). Second order admissible
estimation of variance. {\it Statistics and Decisions, supplement issue \bf 3},
17--29.
}

\Ref{
Malley, J.D. and Hornstein, J. (1993). Quantum statistical inference.
{\it Statistical Science \bf 8}, 433--457.
}

\Ref{
Massar, S. and Popescu, S. (1995).
Optimal extraction of information from finite quantum ensembles.
{\it Physical Review Letters} {\bf 74} 1259--1263.
}

\Ref{
Maudlin, T. (1994).
{\it Quantum Non-locality and Relativity}.
Blackwell, Oxford.
}

\Ref{
Meyer, P.A. (1986). El\'ements de probabilit\'es quantiques.
pp.\ 186--312 in: {\it S'eminaire de Probabilit\'es XX},
ed.\ J. Az\'ema and M. Yor,
Springer Lecture Notes in Mathematics {\bf 1204}.
}

\Ref{
Penrose, R. (1994).
{\it Shadows of the Mind: a Search for the Missing Science of Consciousness}.
Oxford University Press.
}

\Ref{
Percival, I. (1998).
{\it Quantum State Diffusion}.
Cambridge University Press.
}

\Ref{
Peres, A. (1995).
{\it Quantum Theory: Concepts and Methods}.
Kluwer, Dordrecht.
}

\Ref{
Peres, A. and Wootters, W.K. (1991).
Optimal detection of quantum information.
{\it Physical Review Letters \bf 66} 1119--1122.
}

\Ref{
Schipper, C.M.A. (1997).
{\it Sharp Asymptotics in Nonparametric Estimation}.
PhD thesis, University Utrecht, ISBN 90-393-1208-7.
}

\Ref{
van Trees, H.L. (1968).
{\it Detection, Estimation and Modulation Theory (Part 1).}
Wiley, New York.
}

\Ref{
Vidal, G., Latorre, J.I., Pascual, P., and Tarrach, R. (1998).
{\it Optimal minimal measurements of mixed states}.
\newline Preprint {\tt quant-ph/9812068}, {\tt http://xxx.lanl.gov}.
}

\newpage

\address
{Mathematical Institute\\
University Utrecht\\
P.O. Box 80010\\
3508 TA Utrecht\\
Netherlands\\}
\email {gill@math.uu.nl}

\end{document}